\documentstyle{amsppt}
\magnification=\magstep1
\hsize=6 truein
\hcorrection{.375in}
\vsize=8.5 truein
\parindent=20pt
\baselineskip=14pt
\TagsOnRight
\NoBlackBoxes
\footline{\hss\tenrm\folio\hss}

\centerline{On Invariant Measures}
\centerline{for $Diff(S^1)$}

\vskip1truein

\centerline{Doug Pickrell}
\centerline{Mathematics Department}
\centerline{University of Arizona}
\centerline{Tucson, Arizona 85721 }

\vskip1truein

\flushpar Abstract.  In [Pi] we constructed biinvariant measures 
(possibly having values in a line bundle) for a unitary loop 
group $LK$ acting on the formal completion of its 
complexification $LG$. One motivation for this was to find a
geometric construction for the unitary structure of the positive
energy representations of $LK$. In this paper we pursue an analogous 
construction for $Diff(S^1)$. 

\vfill
\eject

\centerline{\S 0. Introduction.}

\bigskip

The complex Virasoro algebra is the universal central extension 
of the Lie algebra of complex trigonometric vector fields on 
the circle.  As a vector space 
$$Vir=(\sum_{n\in \Bbb Z}\Bbb CL_n)\oplus \Bbb C\kappa ,\tag 0.1$$
where $L_n=ie^{in\theta}\frac d{d\theta}=-z^{n+1}\frac d{dz}$; the bracket is determined by 
the relations
$$[L_n,L_m]=(m-n)L_{n+m}+\frac 1{12}n(n^2-1)\delta (n+m)\kappa ,\quad 
[L_n,\kappa ]=0.\tag 0.2$$
The universal central extension $\hat {\Cal D}_{an}$ of $\Cal D_{
an}$, the analytic, 
orientation-preserving diffeomorphisms of $S^1$, is a global real 
form for $Vir$.  

The Virasoro algebra has a triangular decomposition, in the 
technical sense of [MP], where
$$\frak n^{\pm}=\sum_{\pm n>0}\Bbb CL_n,\quad and\quad\hat {\frak h}
=\Bbb CL_0\oplus \Bbb C\kappa ,\tag 0.3$$
so that $Vir$ and $\hat {\Cal D}_{an}$ are, in some limited respects, similar to a 
rank two affine Kac-Moody Lie algebra $\frak g(A)$ and a global real 
form $K(A)$.  In particular the highest weight representations 
can be realized via a Borel-Weil type construction ([KY]), and 
the unitarizable modules can be exponentiated to $\hat {\Cal D}_{
an}$.  For this 
reason, following [Pi], it is natural to inquire whether there 
exist invariant (bundle-valued) measures for $\hat {\Cal D}_{an}$, acting on a 
completion of $\Cal D_{an}$.  This question arises in string theory as 
well, where the appropriate completion is $\Cal D_{qs}$, the group of 
quasisymmetric homeomorphisms of $S^1$ (see the ``physicist's 
wish-list'' in \S 5 of [Pe]).  In this paper we will construct 
(bundle-valued) measures, using basically the same technique as 
in Part III of [Pi].  We conjecture that these measures are 
invariant, and we would like to believe that they are supported 
on $\Cal D_{qs}$, but we are very far from being able to prove this 
rigorously.  

Only fragments of the general setup in [Pi] carry over 
to the Virasoro context.  In particular there is not a complex 
algebraic Virasoro group $\hat {\Cal D}^{\Bbb C}$, so that it is unlikely that there 
exists a good analogue of the formal completion, $\hat {\Cal D}^{
\Bbb C}_{formal}$.  
Despite this we can consider the top stratum
$$\Cal N^{-}\times\hat {\Cal H}\times \Cal N^{+}\subset\hat {\Cal D}^{
\Bbb C}_{formal},\tag 0.4$$
where $\Cal N^{\pm}$ are the simply connected profinite nilpotent Lie 
groups corresponding to the Lie algebras 
$$\frak n^{\pm}_{formal}=\prod_{\pm k>0}\Bbb CL_k,\tag 0.5$$
respectively, and $\hat {\Cal H}$ is the simply connected group 
corresponding to $\hat {\frak h}$.  Unlike the situation for Kac-Moody 
groups, there is an injection 
$$\hat {\Cal D}\to \Cal N^{-}\times\hat {\Cal H}\times \Cal N^{+}
.\tag 0.6$$
In concrete terms this amounts to the classical fact that each 
(quasisymmetric) $\phi\in \Cal D$ can be uniquely factored as a 
composition of maps 
$$\phi =l\circ diag\circ u\tag 0.7$$
where 
$$u=z(1+\sum_{n\ge 1}u_nz^n)\tag 0.8$$
is a univalent holomorphic function on the disk, $diag$ is 
rotation by a complex number of absolute value $\le 1$, the 
mapping inverse to $l$, 
$$l^{-1}=z+\sum_{n\ge 0}b_nz^{-n},\tag 0.9$$
is a univalent holomorphic function on the disk about infinity, 
and the compatibility condition
$$diag(u(S^1))=l^{-1}(S^1)\tag 0.10$$
is satisfied (see for example page 100 of [L]).  Similarly there 
does not exist a formal flag space (in the sense of [Pi]), but 
because of the inclusion (0.6), and the Bieberbach-de Branges 
inequalities ([B]), we can reasonably define the formal 
completion of $Rot(S^1)\backslash \Cal D_{an}$, viewed as a subset of $
\Cal N^{+}$, as 
$$(Rot(S^1)\backslash \Cal D)_{formal}=\{u=z(1+\sum_{n\ge 1}u_nz^
n):\vert u_n\vert <n+1\},$$
where we are viewing $u$ simply as a formal power series. The 
Lie algebra $vir$ acts holomorphically on these spaces in a 
natural way.

As in the theory of loop groups, the elemental matrix 
coefficients can be formally written in terms of Toeplitz 
determinants 
$$\sigma_{c,h}(\hat{\phi })=det(A_a(\hat{\phi }))^{c-8h}det(A_p(\hat{
\phi }))^{8h},\tag 0.11$$
$$=det(A_a(\hat{\phi }))^cdiag(\phi )^{8h}$$
where $\phi$ denotes a diffeomorphism of $S^1$, $\hat{\phi}$ is in the universal 
central extension, and the subscripts refer to a choice of periodic 
or antiperiodic spin structure for $S^1$.  To make sense of the 
formal measure 
$$\vert diag(\phi )\vert^{16h}det\vert A_a(\phi )\vert^{2c}\Cal D
(\phi ),\tag 0.12$$
we proceed as in the case of loop groups, by first constructing 
a regularization with parameter 
$$d\nu_{\beta ,c,h}(\phi )=\frac 1E\vert diag(\phi )\vert^{16h}de
t\vert A_a(\phi )\vert^{2c}d\nu_{\beta}(\phi ),\tag 0.13$$
where $\nu_{\beta}$ denotes a (slight variation of a) left quasi-invariant 
measure on $\Cal D_{C^1}$ constructed by the Malliavins in [MM] and 
Shavgulidze in [Sh].  One can make sense of this regularization 
by essentially the same method used in the loop group case 
(for $c,h\ge 0$).  From a geometric point of view, one should 
think of $\sigma_{c,h}$ as a section of a line bundle $\Cal L^{*}_{
c,h}\to \Cal D$; the 
expression (0.12) is heuristically 
$$\langle\sigma_{c,h},\sigma_{c,h}\rangle \Cal D\phi ,\tag 0.14$$
the density for the norm squared of $\sigma_{c,h}$ as a section.  From 
this point of view, $\nu_{\beta ;c,h}$ is a coordinate expression for a 
measure having values in the line bundle 
$$\vert \Cal L_{c,h}\vert^2=\vert DetA_a\vert^{-2(c-8h)}\otimes\vert 
DetA_p\vert^{-16h}\to \Cal D\tag 0.15$$
It is very likely that the bundle-valued measure $\nu_{\beta ;c,h}$ is 
asymptotically invariant as $\beta\downarrow 0$, with respect to $
\Cal D_{an}$ acting 
from the left.  This is an interesting issue, because 
it turns out to be intimately related to stationary phase 
approximations for Brownian paths.  The measure $\nu_{\beta ;c,h}$ is not 
quasi-invariant for the right action, although it does 
heuristically appear to be invariant in the limit as $\beta\downarrow 
0$.  

The next step is to take the limit of $\nu_{\beta ,c,h}$ as the inverse 
temperature $\beta$ tends to zero.  The Bieberbach-de Branges 
inequalities (or just inequalities coming from the area theorem) 
immediately imply that weak limits exist.  Unfortunately in 
the Virasoro case there is not enough structure to easily parlay this 
into a statement about the existence of an invariant measure $\mu$ 
with values in $\vert \Cal L_{c,h}\vert^2$.  What we lack, in contrast to the 
Kac-Moody case, is a priori knowledge of certain elemental 
distributions (the first nontrivial coefficients $u_1$ and $b_0$ for the 
univalent functions $u$ and $l$ above, which correspond to the 
simple roots). 

If we do succeed in proving that the limit
$$\mu_{c,h}=\lim_{\beta\downarrow 0}\nu_{\beta ,c,h}\tag 0.16$$
exists and is $\Cal D_{an}$-invariant (as a bundle-valued measure), then 
we will have an appealing geometric way of understanding why 
the lowest weight module corresponding to $(c,h)$ is unitary.  It 
is intuitively clear that this should work for only a continuous 
set of parameters, and unfortunately our work so far does not 
pinpoint where the constraint $c>1$ should arise (we suspect 
that it arises in the calculation of the distribution for the $u_
1$ 
coefficient - this probably has an absolutely continuous 
distribution only for $c>1$). 

\bigskip

\centerline{\S 1. Completions and Classical Analysis.}

\bigskip

Virtually all the material in this section is well-known.
In constructing the loop group analogue of the measure (0.7), 
Peller's work on the Schatten properties of Hankel operators 
for multiplication operators played a pivotal role.  The 
analogous results for homeomorphisms of $S^1$ have been 
established by Coifman, Meyer, Peller, Rochberg and Semmes.  
In \S 1.2 we formulate their fundamental work in the context of 
this paper.  In \S 1.5 we discuss embeddings into the top stratum 
of the (nonexistent) formal completion, which involves the
Ahlfors-Beurling theory of quasi-symmetric homeomorphisms of 
$S^1$.

\bigskip

\flushpar\S 1.1. Notation.

\smallskip

We will denote a group of orientation preserving 
homeomorphisms of $S^1$ by $\Cal D$, where a subscript will indicate 
the order of smoothness of the elements of the group.  The 
universal covering will be denoted by $\tilde {\Cal D}$, and it will be 
identified with the set of $\Phi :\Bbb R\to \Bbb R$ such that 
$$\Phi (s+2\pi )=\Phi (s)+2\pi ,\quad\Phi\quad is\quad increasing
,\tag 1.1.1$$
and $\Phi$ $\underline {and}$ $\Phi^{-1}$ satisfy the appropriate smoothness condition.  
We will frequently write a homeomorphism of $S^1$ as 
$$\phi (e^{i\theta})=e^{i\Phi (\theta )},\tag 1.1.2$$
without further comment.

The group $\Cal D_{C^{\infty}}$ acts naturally and unitarily on half-densities 
on $S^1$.  We identify these densities with functions in the usual 
way, 
$$L^2(S^1)\to\Omega^{1/2}_{L^2}(S^1):f\to f\vert d\theta\vert^{1/
2}.\tag 1.1.3$$
Then $\phi\in \Cal D_{C^{\infty}}$ corresponds to a unitary operator $
U_{\phi}$ of $L^2(S^1)$, 
where 
$$U_{\phi^{-1}}:f\to (\Phi')^{1/2}f\circ\phi .\tag 1.1.4$$

Relative to the Hardy space polarization of $L^2(S^1)$, we will write
$$U_{\phi}=\left(\matrix A&B\\C&D\endmatrix \right)\tag 1.1.6$$

In the next subsection we will be interested in intrinsic 
characterizations of certain induced topologies on $\Cal D_{C^{\infty}}$.  A 
standard prototype is the following standard

\proclaim{(1.1.7)Proposition}The closure of $\Cal D_{C^{\infty}}$ in the strong 
operator topology for $U(L^2(S^1))$ is
$$\Cal D_{W^1_{L^1}}=\{\phi\in \Cal D:\phi\quad and\quad\phi^{-1}
\quad are\quad absolutely\quad continuous\}.$$
Thus $U_{\phi_n}\to U_{\phi}$ strongly if and only if $\phi_n\to\phi$ and $
\psi_n\to\psi$ in 
$W^1_{L^1}$, where $\psi =\phi^{-1}$.
\endproclaim

\flushpar(1.1.8)Remarks. (a) It is not the case that $\phi$ absolutely 
continuous implies that its inverse $\psi$ is absolutely continuous. 
For example, let $S$ denote a generalized Cantor set of positive 
measure. Then the complement $S^c$ is open and dense, hence
$$\Phi (x)=\int_0^x\chi_{S^c}ds$$
is strictly increasing and absolutely continuous.  But $S$ has 
positive measure while $\Phi (S)$ has measure zero, hence its 
inverse is not absolutely continuous. This illustrates why we 
generally need to impose conditions on both $\phi$ and its inverse.
In a vague way, this explains why the conditions that arise 
in the next subsection are expressed in terms of $log\Phi'$.

(b) A second, but less interesting, topology on $\Cal D_{W^1_{L^1}}$ is that 
induced by the operator norm topology on $U(L^2)$. In this case 
the induced topology is the discrete topology.

\bigskip

\flushpar$\S 1.2$.  The $p$-Kac-Peterson completion and restricted 
groups.  

\smallskip

Given a polarized separable Hilbert space, $H=H_{+}\oplus H_{-}$, and a 
symmetrically normed ideal $\Cal I$, there is an associated Banach 
$*$-algebra, $\Cal L_{(\Cal I)}$, which consists of bounded operators on $
H$, 
represented as two by two matrices with respect to the 
polarization, as in (1.1.5), such that the norm 
$$\vert\left(\matrix A&\\&D\endmatrix \right)\vert_{\Cal L}+\vert\left
(\matrix &B\\C&\endmatrix \right)\vert_{\Cal I}\tag 1.2.1$$
is finite.  The $*$-operation is the usual adjoint operation.  The 
corresponding unitary group is 
$$U_{(\Cal I)}=U(H)\cap \Cal L_{(\Cal I)};$$
it is referred to as a restricted unitary group in [PS].  
Geometrically this group is the group of automorphisms of a 
Grassmannian (Finsler) symmetric space modelled on $\Cal I$.  There 
are two obvious topologies on $U_{(\Cal I)}$.  The first is the induced 
Banach topology, and in this topology $U_{(\Cal I)}$ has the additional 
structure of a Banach Lie group.  The second is the Polish 
topology $\tau_{KM}$ for which convergence means that for 
$g_n,g\in U_{(\Cal I)}$, $g_n\to g$ if and only if $g_n\to g$ strongly and 
$$\left(\matrix &B_n\\C_n&\endmatrix \right)\to\left(\matrix &B\\
C&\endmatrix \right)\quad in\quad \Cal I.$$
For $\Cal I=\Cal L_2$, Hilbert-Schmidt operators, the identity component 
of $(U_{(\Cal I)},\tau_{KM})$ is the Kac-Peterson completion of the (infinite 
classical) Kac-Moody group of type $A_{2\infty}$, modulo its center (see 
the introduction to Part II of [Pi]).  

When $H=L^2(S^1,\Bbb C^n)$ with the Hardy polarization, and $\Cal I
=\Cal L_p$, 
there is a beautiful intrinsic description of the (metric 
isomorphism class of the) induced norm on the subalgebra 
$$Map(S^1,\Cal L(\Bbb C^n))\subset \Cal L_{(\Cal L_p)};\tag 1.2.2$$
in this case Peller has shown that the norm (1.2.1) is 
equivalent to 
$$\vert F\vert_{L^{\infty}}+\vert F\vert_{B^{1/p}},\tag 1.2.3$$
where the Besov p-norm is defined by (3.3.3) of Part III (if 
$p=\infty$, then the Besov space is replaced by $VMO$, and the 
result in this case is due to Hartmann; to my knowledge it not 
known how to intrinsically characterize the induced norm 
corresponding to other ideals of interest, such as the Dixmier 
trace class).  This implies that the map of Banach Lie groups 
$$L_{L^{\infty}\cap B^{1/p}}U(n)\to U_{(\Cal L_p)}\tag 1.2.4$$
is a homeomorphism onto its image (in finite dimensions, it
would automatically follow that this map is an embedding, but 
this does not appear to be true here, because for the map of 
Lie algebras, the image does not appear to have a complement). 
For the Kac-Peterson type topology,
$$L_{B^{1/p}}U(n)\to (U_{(\Cal L_p)},\tau_{KM})$$
is a homeomorphism onto its image (when $p=\infty$, we replace 
the Besov class by $VMO$), and when $p=2$ this describes the 
Kac-Peterson completion of the loop group.  

To formulate the analogue of this result for $\Cal D$, we will have to 
work exclusively with the Kac-Moody type topologies, because 
the inclusion 
$$\Cal D_{C^{\infty}}\to U_{(\Cal L_p)}\tag 1.2.5$$
defined by (1.1.2-5) is not continuous for the Banach topology 
(as in (b) of (1.1.8), the induced topology is the discrete topology).  
From now on it will be understood that $U_{(\Cal L_p)}$ is equipped with 
the topology $\tau_{KM}$.  We will refer to the completion of the 
smooth diffeomorphism group in $U_{(\Cal L_2)}$ as the Kac-Peterson 
completion, since this is the analogue of the corresponding 
object for Kac-Moody groups.  The problem then is to 
intrinsically characterize the Kac-Peterson completion and its 
$p$-generalization.  

The following fundamental theorem (with $\Bbb R$ in place of $S^1$) 
first appeared in full generality in the dissertation of Semmes 
(Theorem 3 of [S]).  

\proclaim{(1.2.6)Theorem} (a) For $p=\infty$,
$$\Cal D_{W^1_{L^1}}\cap U_{(\Cal L_{\infty})}=\{\phi\in \Cal D:l
og\Phi'\in VMO\}.$$

\flushpar(b) For $1\le p<\infty$,
$$\Cal D_{W^1_{L^1}}\cap U_{(\Cal L_p)}=\{\phi\in \Cal D:log\Phi'
\in B^{1/p}\}.$$
\endproclaim

The direction of critical importance to us is containment.
The proof proceeds as follows.  Suppose that $log\Phi'\in B^{1/p}$ (or 
$VMO$ if $p=\infty$).  We must show that 
$$U_{\phi}\circ j\circ U_{\phi}^{-1}-j\in \Cal L_p.\tag 1.2.7$$
Coifman and Meyer considered the homotopy 
$$\Phi_t(\theta )=\Phi (0)+\int_0^{\theta}\Phi'(\tau )^td\tau ,$$
from a rigid rotation to $\Phi$.  They differentiated the 
corresponding deformation of (1.2.7) to obtain 
$$U_{\phi}\circ j\circ U_{\phi}^{-1}-j=\int_0^1(U_{\phi_t}\circ L_
2(log\Phi'\circ\Phi_t^{-1})\circ U_{\phi_t}^{-1})dt\tag 1.2.8$$
where $log\Phi'\circ\Phi_t^{-1}\in B^{1/p}$, $\forall t$, and
$$L_2(b)=[B\frac d{d\theta},j]+\frac 12[b,j],\quad B'=b.\tag 1.2.9$$
Rochberg and Peller (Coifman and Meyer, and Hartmann, resp.)  
obtained the necessary $\Cal L_p$ estimates on the first and second 
terms of $L_2$, respectively, for $p<\infty$ (for $p=\infty$, resp.).  

\bigskip

\flushpar\S 1.3. Spin structures.

\smallskip

The circle $S^1$ has two distinct real spin structures, periodic 
(or trivial) and antiperiodic (or Mobius).  In the latter case 
there is a natural action by $\Cal D_{an}^{(2)}$, the double cover.  The 
complexification of the antiperiodic spin structure is trivial, 
but not equivariantly trivial.  In each case there is a natural 
Hilbert space structure for half-forms, denoted by $H_p$ and $H_a$, 
respectively.  We will identify both of these spaces with $H$ as 
follows.  In the periodic case the identification is simply 
$$H\to H_p:f\to f(d\theta )^{1/2}.\tag 1.3.1$$
In the antiperiodic case there is a polarization 
$$H_a=H_a^{+}\oplus H_a^{-},\tag 1.3.2$$
where $H_a^{\pm}$ is the closure of holomorphic sections of the spin 
bundle for the disk $\bar {D}^{\pm}$, respectively. There is an isomorphism 
of polarized spaces
$$H\to H_a:f\to f(dz)^{1/2}.\tag 1.3.3$$
We let
$$\Cal D_{W^1_{L^1}}^{(2)}\to U(H):\tilde{\phi}\to V_{\tilde{\phi}}
=\left(\matrix A_a(\tilde{\phi })&B_a(\tilde{\phi })\\C_a(\tilde{
\phi })&D_a(\tilde{\phi })\endmatrix \right)\tag 1.3.4$$
denote the induced action on $H$.  The same chain of arguments 
as in the previous subsection (one simply systematically 
replaces $d\theta$ by $dz$) implies the following 

\proclaim{(1.3.5)Theorem}For each $p\ge 1$,  
$$V_{\tilde{\phi}}\in U_{(p)}\Leftrightarrow U_{\phi}\in U_{(p)}\Leftrightarrow 
log\Phi'\in B^{1/p}.$$
\endproclaim

\flushpar(1.3.6)Remark. In the next section we will see that
$$det\vert A_p(\phi )\vert\le det\vert A_a(\tilde{\phi })\vert$$
for all $\phi$. It is natural ask whether this actually follows from 
an inequality of the operators, i.e. is
$$\vert A_p(\phi )\vert\le\vert A_a(\tilde{\phi })\vert\quad ?$$

\smallskip

For later purposes it is convenient to briefly discuss the 
implication 
$$log\Phi'\in B^{1/2}\Rightarrow V_{\tilde{\phi}}\in U_{(2)},\tag 1.3.7$$
by considering the square-integrability of the kernel for $[j,V_{
\phi}]$.  
Using the fact that 
$$proj_{H_{\pm}}:H\to H_{\pm}:f\to\frac {\pm 1}{2\pi i}\int_{S^1}
f(\zeta )\frac 1{\zeta -z}d\zeta ,\quad\vert z\vert <1\quad (resp
.\quad\vert z\vert >1),$$
a direct computation shows that the kernel for $\frac 12[j,V_{\phi}
]$ is 
given by 
$$-\frac 1{2\pi i}(\frac {^1}{z-\phi (\zeta )}(\frac {d\phi}{d\zeta}
)^{1/2}-\frac 1{\psi (z)-\zeta}(\frac {d\psi}{dz})^{1/2})d\zeta ,\tag 1.3.8$$
where $\psi =\phi^{-1}$.  To check that this kernel is $L^2$, it is 
equivalent and slightly easier to check that the kernel for 
$V_{\phi^{-1}}\circ\frac 12[j,V_{\phi}]$ is $L^2$, and this is given by 
$$K_{\phi}(s,t)\frac {d\zeta}{2\pi i}=-\frac 1{2\pi i}(\frac {(\frac {
d\phi}{d\zeta}\frac {d\psi}{dz})^{1/2}(z-\zeta )-(\phi (z)-\phi (
\zeta ))}{(\phi (z)-\phi (\zeta ))(z-\zeta )})d\zeta$$
$$=-\frac 1{2\pi i}(\frac {(\frac {\phi (z)}z\Phi'(t)\frac {\phi 
(\zeta )}{\zeta}\Phi'(s))^{1/2}\frac {z-\zeta}{\phi (z)-\phi (\zeta 
)}-1}{z-\zeta})d\zeta$$
$$=\frac 1{2\pi i}\frac {(\Phi'(s)\Phi'(t))^{1/2}\frac {sin((t-s)
/2)}{sin((\Phi (t)-\Phi (s))/2)}-1}{z-\zeta}d\zeta .$$
$$=\frac 1{2\pi i}\frac {exp(\frac {B(t)+B(s)}2-ln(\frac {sin(\frac 
12\int_s^te^B)}{sin((t-s)/2)}))-1}{z-\zeta}d\zeta\tag 1.3.9$$
where $z=e^{it},\zeta =e^{is}$, and $B(t)=ln\Phi'(t)$.  The proof of (3.3.6) 
below implicitly indicates how to directly show this is an $L^2$ 
kernel, assuming that $B\in B^{1/2}$.  It would be interesting to 
know if this follows from a general result.  

\bigskip

\flushpar\S 1.4. The symplectic action.

\smallskip

Let $V$ denote the symplectic vector space
$$V=(Map_{W^{1/2}}(S^1,\Bbb R)/\Bbb R,\omega ),\quad where\quad\omega 
(f,g)=\int fdg.\tag 1.4.1$$
In analogy with (1.1.7), Nag and Sullivan ([NS]) have proven
the following

\proclaim{(1.4.2)Theorem}With respect to the natural action of 
$\Cal D$ on measureable functions $(\Cal C_{\phi}:f\to f\circ\phi^{
-1})$, we have 
$$\Cal D\cap Sp(V)=\Cal D_{qs},$$
the group of quasi-symmetric homeomorphisms.
\endproclaim

An important consequence of this is the following.  With 
respect to the Hardy polarization of $V^{\Bbb C}$, where $V^{\Bbb C}_{
\pm}$ are the 
conformally invariant Hilbert spaces 
$$V^{\Bbb C}_{\pm}=\{f\in H^0(\{\pm\vert z\vert <1\}):\frac i2\int 
df\wedge d\bar {f}<\infty \}/\Bbb C\tag 1.4.3$$
we can write the composition operator $\Cal C_{\phi}$ as a two by two 
matrix as in (1.1.6).  Since $\Cal C_{\phi}$ is symplectic, $A(\Cal C_{
\phi})$ is always
invertible.

Since the Hardy polarization is positive, we can form the 
restricted symplectic group 
$$Sp_{(\Cal I)}=\{\left(\matrix A&B\\C&D\endmatrix \right)\in Sp(
V):B,C\in \Cal I\},\tag 1.4.4$$
for any symmetrically normed ideal (see [Se1]).  A corollary of 
(1.3.5) is that 
$$\Cal D\cap Sp_{(\Cal L_p)}=\{\phi :ln\Phi'\in L^{\infty}\cap B^{
1/p}\}.\tag 1.4.5$$

\bigskip

\flushpar\S 1.5. Mapping $\Cal D$ into the Formal Completion $\Cal D_{
formal}^{\Bbb C}$.

\smallskip

The formal completion $\Cal D^{\Bbb C}_{formal}$, if it did exist, would contain 
as a dense open subset the product space 
$$\Cal N^{-}\cdot \Bbb C^{*}\cdot \Cal N^{+}\subset \Cal D^{\Bbb C}_{
formal},\tag 1.5.1$$
where $\Cal N^{\pm}$ are the simply connected profinite nilpotent groups 
corresponding to $\frak n_{formal}^{\pm}$, respectively, and $\Bbb C^{
*}$ is identified 
with the complexification of $Rot(S^1)$.  In this subsection, after 
some preliminary discussion, we will show that various 
analytical completions of $\Cal D$ (and also the Neretin-Segal 
semigroup of annuli) can be naturally mapped into $\Cal N^{-}\cdot 
\Bbb C^{*}\cdot \Cal N^{+}$.
 
First we must realize the group $\Cal N^{+}=exp(\frak n^{+}_{form
al})$ in a 
concrete way.  There is a natural Lie algebra representation of 
$\frak n^{+}_{formal}$ by continuous derivations of the algebra $
\Bbb C[[z]]$, 
$$\frak n^{+}_{formal}\times \Bbb C[[z]]\to \Bbb C[[z]]:(V,f)\to 
v\cdot f',\tag 1.5.2$$
where $V=v\frac {\partial}{\partial z}$, $v=O(z^2)$ (continuous means 
$\sigma (\sum f_nz^n)=\sum f_n\sigma (z^n)$).   We claim that this can be 
exponentiated to a group action of $\Cal N^{+}$ by continuous algebra 
automorphisms 
$$\Cal N^{+}\times \Bbb C[[z]]\to \Bbb C[[z]].\tag 1.5.3$$
To see this, we compute that
$$e^Vf=\sum_{n,m\ge 0}\frac 1{n!}f_m(v\frac {\partial}{\partial z}
)^n(z^m).\tag 1.5.4$$
By induction it is easy to see that
$$(v\frac {\partial}{\partial z})^n=\sum_{1\le l\le n}w_{n,l}(\frac {
\partial}{\partial z})^l,\quad where\quad w_{n,l}=O(z^{n+l}).\tag 1.5.5$$
This implies that
$$(v\frac {\partial}{\partial z})^n(z^m)=O(z^{n+l+m-l})=O(z^{n+m}
),\tag 1.5.6$$
and hence the right hand side of (1.5.4) is a well-defined formal 
power series.  Given that $e^V$ is well-defined, it is routine to 
check that $e^V$ is a continuous automorphism, and 
$$\Cal N^{+}\to Aut_{C^0}(\Bbb C[[z]])\tag 1.5.7$$
is a faithful representation.

The image of (1.5.7) can be described geometrically in the following 
way.  

\proclaim{(1.5.8)Proposition}With reference to the group action 
(1.5.3), 
$$\Cal B^{+}=\Bbb C^{*}\cdot \Cal N^{+}=Aut_{C^0}(\Bbb C[[z]])=(\{
u\in \Bbb C[[z]]:u=\lambda z(1+\sum_{n\ge 1}u_nz^n),\lambda\ne 0\}
,\circ ),$$
where the multiplication in the latter group is composition of 
formal power series, and $u$ acts by
$$u:\Bbb C[[z]]\to \Bbb C[[z]]:f\to f\circ u^{-1}.$$
\endproclaim

\flushpar Proof of (1.5.8).  Given a continuous automorphism $\sigma$, 
we must produce $\lambda\in \Bbb C^{*}$, $V=v\frac {\partial}{\partial 
z}\in \frak n^{+}_{formal}$ and $u$ such that 
$$\lambda e^V=\sigma =(\cdot )\circ u^{-1}.\tag 1.5.9$$
Since $\sigma$ is continuous, it is determined by what it does to $
z$, so 
that clearly $u$ is determined by $\sigma (z)=u^{-1}$. To produce $
\lambda$ and $V$,
note that for $n>1$, in the expansion (1.5.5), 
$$w_{n,1}=(v\frac {\partial}{\partial z})^{n-1}(v).\tag 1.5.10$$
Applying $\lambda e^V$ to $z$, we obtain the equation
$$\lambda (z+\sum_{n>0}\frac 1{n!}(v\frac {\partial}{\partial z})^{
n-1}(v))=u^{-1}=c_1z+c_2z^2+....\tag 1.5.11$$
This is equivalent to relations
$$\lambda =c_1,\quad\lambda v_2=c_2,\quad\lambda (v_3+v_2^2)=c_3,
\quad\lambda (v_4+\frac 52v_3v_2+\frac 13v_3^3)=c_4,\tag 1.5.12$$
and so on, from which it is clear that we can obtain $V$ 
(although I do not see how to write these relations down in 
simple closed form).//

\smallskip

\flushpar(1.5.13) Remark. There is a similar identification
$$\Cal B^{-}=\Bbb C^{*}\cdot \Cal N^{-}=Aut_{C^0}(\Bbb C[[z^{-1}]
]).$$
\smallskip

Now we consider $\Cal D$, and more generally the Neretin-Segal 
semigroup of annuli $\Cal A$.  The appropriate smoothness class is 
the Ahlfors-Beurling group of quasi-symmetric 
orientation-preserving homeomorphisms of $S^1$, $\Cal D_{qs}$.   The 
corresponding smoothness condition for Segal's category of 
compact Riemann surfaces with parameterized boundaries is the 
following.  An object $\Sigma$ should be an oriented compact $C^0$ 
surface with boundary, with a compatible complex structure in 
the interior, and a continuous parameterization of each of its 
boundary components.  In addition for each point $q\in\partial\Sigma$, there 
should be a chart 
$$z:U\to \{\vert z\vert\le 1,\Bbb Re(z)\ge 0\}$$
which is holomorphic in the interior, such that each 
parameterization is quasi-symmetric when expressed in any 
pair of such charts. In this context sewing of such objects is 
well-defined.

Now consider the semigroup of annuli $\Cal A_{qs}$.  Given $A\in 
\Cal A_{qs}$, 
there is a unique holomorphic isomorphism 
$$\bar {D}^{-}\circ A\circ\bar {D}^{+}=\hat {\Bbb C}\tag 1.5.14$$
which is characterized by the condition that restriction 
induces maps satisfying
$$\lambda u:\bar {D}^{+}\to\hat {\Bbb C},\quad u(z)=z(1+\sum_{n\ge 
1}u_nz^n)\tag 1.5.15$$
$$l^{-1}:\bar {D}^{-}\to\hat {\Bbb C},\quad l^{-1}(z)=z+\sum_{n\ge 
0}b_nz^{-n},\tag 1.5.16$$
where $\lambda\in \Bbb C^{*}$.

\proclaim{(1.5.17)Proposition} There is an equivariant inclusion
$$\Cal A_{qs}\to \Cal N^{-}\cdot \Bbb C^{*}\cdot \Cal N^{+}:A\to 
(l,diag,u),$$
where $l$ is the mapping inverse to $l^{-1}$, and $diag=Rot(\lambda 
)$.
By passing to the portion of the boundary where 
$\lambda u(S^1)=l^{-1}(S^1)$, this also induces an equivariant inclusion 
$$\Cal D_{qs}\to \Cal N^{-}\cdot \Bbb C^{*}\cdot \Cal N^{+}:\phi\to 
(l,diag,u),$$
where as a mapping $\phi (z)=l\circ diag\circ u(z)$, $z\in S^1$, and $
l$ and $u$ 
have quasi-conformal extensions to $\hat {\Bbb C}$.  This factorization of $
\phi$ 
is unique.  \endproclaim

\flushpar(1.5.18)Remarks.  (a) The meaning of equivariant is the 
following.  Although $\Cal D_{an}$ does not act on $\Cal N^{-}\cdot 
\Bbb C^{*}\cdot \Cal N^{+}$, the Lie 
algebra of polynomial vector fields on $S^1$ does act by 
holomorphic vector fields from both the left and the right.  
The formulas are exactly the same as in the Kac-Moody case.  
These are written down in \S 1.6 of Part I of [Pi] (see especially 
(1.6.6)).  Similarly there is an action of $vect_{poly}$ on the space 
of annuli (see \S 2 of [Se2]). These actions commute with the 
maps in (1.5.17).

(b) The uniqueness assertion can be understood in the following 
way.  Suppose that $\phi\in \Cal D_{qs}$.  Then we know that $\phi
\in Sp(V)$, as 
in \S 1.4, and hence $A(\Cal C_{\phi})$ is an invertible operator.  Since $
l$ must 
have the form $l=z+\sum_{n\ge 0}l_nz^{-n}$, it follows that we must have 
$$diag\cdot u=A(\Cal C_{\phi})^{-1}(z).\tag 1.5.19$$

(c) Note that this argument also proves a weak version of the 
existence of the factorization.  For given (1.5.19), $u$ has $W^{
1/2}$ 
boundary values, and we also obtain $l^{-1}$, which also has $W^{
1/2}$ 
boundary values, and on $S^1$ 
$$l^{-1}=diag\cdot u\circ\phi^{-1}.\tag 1.5.20$$
To see that $diag\cdot u$ is univalent, suppose that $\vert z_0\vert 
<1$.  Then 
because $u(z_0)$ is not in $\frac 1{diag}l^{-1}(\hat {\Bbb C}\setminus 
D)$, 
$$\frac 1{2\pi i}\int_{S^1}dlog(u-u(z_0))=\frac 1{2\pi i}\int_{S^
1}dlog(\frac 1{diag}l^{-1}\circ\phi -u(z_0))$$
$$=\frac 1{2\pi i}\int_{S^1}dlog(\frac 1{diag}l^{-1}-u(z_0))=1.$$
Thus $u$ is univalent, and similarly $l^{-1}$ is univalent.  
Unfortunately this simple argument does not explain why $u$ and 
$l^{-1}$ are continuous on $S^1$, let alone that they have 
quasi-conformal extensions.  

(d) The area of the annulus $A$ is given by
$$area(A)=\frac i2\int_{S^1}l^{-1}d\bar {l}^{-1}-\frac i2\int_{S^
1}\vert diag\vert^2ud\bar {u}=\pi (1-\sum_{n\ge 0}n(\vert b_n\vert^
2+\vert diag\vert^2\vert u_n\vert^2))\tag 1.5.22$$
where $u_0=1$, so that 
$$\vert b_n\vert^2+\vert diag\cdot u_n\vert^2\le 1/n,\tag 1.5.23$$
and in particular $\vert diag\vert\le 1$.

\bigskip

\centerline{\S 2. The extension $\hat {\Cal D}$ and determinant formulas.}

\bigskip

\flushpar\S 2.1. The universal central extension, $\hat {\Cal D}$.  

\smallskip

The group $\Cal D_{an}$ has a universal central extension 
$$0\to \Bbb Z\oplus i\Bbb R\to\hat {\Cal D}_{an}\to \Cal D_{an}\to 
0.\tag 2.1.1$$
As is well-known ([S1],[KY]), the group $\hat {\Cal D}_{an}$ can be realized in the 
following explicit way.  As a manifold
$$\hat {\Cal D}_{an}=\tilde {\Cal D}_{an}\times i\Bbb R.\tag 2.1.2$$
In these coordinates the multiplication is given by
$$(\Phi ,it)\cdot (\Psi ,is)=(\Phi\circ\Psi ,it+is+iC(\phi ,\psi 
)),\tag 2.1.3$$
where $C$ is the $\Bbb R$-valued cocycle given by
$$C(\phi ,\psi )=\frac 1{48\pi}\Bbb Re\int_{S^1}log(\frac {\partial
\phi}{\partial z}\circ\psi )d(log(\frac {\partial\psi}{\partial z}
))\tag 2.1.4$$
The corresponding Lie algebra is the real form of $vir$ which as 
a vector space equals
$$vect(S^1)\oplus i\Bbb R,\tag 2.1.5$$
with the bracket given by (0.2).

\smallskip

\flushpar(2.1.6)Remarks. (a) Note that the natural domain of $C$ is the 
Kac-Peterson completion, $\{\Phi :ln\Phi'\in B^{1/2}\}$. 

\flushpar(b) To check the assertion about the Lie algebra, 
recall that one obtains the corresponding Lie algebra cocycle 
via
$$c(\vec{\xi },\vec{\eta })=\frac {\partial}{\partial s\partial t}
\vert_{s=t=0}(C(e^{s\vec{\xi}},e^{t\vec{\eta}})-C(e^{t\vec{\eta}}
,e^{s\vec{\xi}}))$$
$$=\frac i{24\pi}\int_{S^1}\frac {\partial\xi}{\partial z}d(\frac {
\partial\eta}{\partial z})=\frac i{24\pi}\int_0^{2\pi}(\tilde{\xi}^{
\prime\prime\prime}(\theta )+\tilde{\xi}'(\theta ))\tilde{\eta }(
\theta )d\theta ,\tag 2.1.7$$
where $\vec{\xi }=\xi (z)\frac d{dz}=\tilde{\xi }(\theta )\frac d{
d\theta}$.  This gives the commutation 
relations in (0.2).  

\smallskip

The homomorphism 
$$\Cal D^{(2)}_{an}\to U_{(2)}:\tilde{\phi}\to V_{\tilde{\phi}},\tag 2.1.8$$
induces a unique homomorphism of central extensions,
$$\hat {\Cal D}_{an}\to\hat {U}_{(2)},\tag 2.1.9$$
where, as in chapter 6 of [PS],
$$\hat {U}_{(2)}=\{(g,q)\in U_{(2)}\times U(H_a^{+}):A_a(g)-q\in 
\Cal L_1\}/U(H_a^{+})_1.\tag 2.1.10$$

It is not known how to write down this homomorphism 
explicitly.  However, it is known that if this homomorphism is 
composed with the basic (lowest weight) representation of $\hat {
U}_{(2)}$ 
(i.e.  the representation described in chapter 10 of [PS]), then 
the module generated by the vacuum is the lowest weight 
module for $vir$ corresponding to the parameters $c=1$, $h=0$.  
Similarly, if we use the homomorphism $\phi\to U_{\phi}$, then we obtain 
the module corresponding to $c=1$ and $h=1/8$ (this is 
calculated for example in section 7 of [S2]).   It follows that 
the matrix coefficient corresponding to the vacuum in the 
representation corresponding to $(c,h)$ is given by (0.11).  

There is an inclusion
$$\hat {\Cal D}_{an}\to \Cal N^{-}\cdot\hat {\Cal H}\cdot \Cal N^{
+}:\hat{\phi}\to l\cdot diag^{\hat{}}\cdot u,\tag 2.1.11$$
where (as a consequence of (0.11))
$$diag^{\hat{}}=(detA_a(\hat{\phi }))^{\kappa}(\frac {detA_a(\hat{
\phi })}{detA_p(\hat{\phi })})^{8L_0},\tag 2.1.12$$
and we are interpreting $detA_a$ and $detA_p$ as the pullbacks of 
the canonical section; this covers the factorization of $\phi$ as a 
composition of functions 
$$\phi =l\circ diag\circ u,\quad diag=(\frac {detA_a(\hat{\phi })}{
detA_p(\hat{\phi })})^{8L_0}=Rot(\frac {detA_p(\hat{\phi })}{detA_
a(\hat{\phi })})^8\tag 2.1.13$$
$u\in \Cal N^{+}$, $l\in \Cal N^{-}$, and $\hat {\Cal H}$ is the simply connected group 
corresponding to $\hat {\frak h}=\Bbb CL_0\oplus \Bbb C\kappa$.  Note that by (1.5.23) we always 
have 
$$\vert diag\vert\le 1,\quad hence\quad det\vert A_p(\phi )\vert\le 
det\vert A_a(\phi )\vert .\tag 2.1.14$$

These formulas have interesting analytical consequences.  For 
each $n>0$, there is an embedding 
$$di_n:sl(2,\Bbb C)\to vir:f\to f_n=\frac {-1}nL_{-n},h\to h_n=\frac 
2nL_0+\frac 1{12}n(n^2-1)\kappa ,e\to e_n=\frac 1nL_n.\tag 2.1.15$$
Geometrically this corresponds to the following. The group of 
projective transformations of $\hat {\Bbb C}$ which map the circle to itself 
is the subgroup $PSU(1,1)\subset PSL(2,\Bbb C)$, where 
$$\left(\matrix a&b\\\bar {b}&\bar {a}\endmatrix \right)\cdot z'=\frac {
\bar {b}+\bar {a}z'}{a+bz'}\tag 2.1.16$$
For $n\ge 1$ there is an $n$-fold covering map, 
$$S^1\to S^1:z\to z'=z^n,\tag 2.1.17$$
and the diffeomorphisms of $z$ which cover the projective 
transformations of $z'$ form a group $PSU(1,1)^{(n)}$ which is an 
$n$-fold covering 
$$0\to \Bbb Z_n\to PSU(1,1)^{(n)}\to PSU(1,1)\to 0\tag 2.1.18$$
The map $di_n$, modulo the center, is the complexification of the 
differential of the embedding 
$$i_n:PSU(1,1)^{(n)}\to \Cal D.\tag 2.1.19$$

Suppose that $\phi\in PSU(1,1)^{(n)}\subset\hat {\Cal D}$ covers $
i_1$ $\left(\matrix a&b\\\bar {b}&\bar {a}\endmatrix \right)\in P
SU(1,1)\subset \Cal D$.  
Then there is a unique factorization 
$$\phi =e^{\frac {-\bar {b}a^{-1}}nL_{-n}}a^{\frac 2nL_0+\frac 1{
12}n(n^2-1)\kappa}e^{\frac {a^{-1}b}nL_{+n}}\in \Cal N^{-}\cdot\hat {
\Cal H}\cdot \Cal N^{+}\tag 2.1.20$$
When we project $\phi$ into $\Cal D$, we obtain a factorization of $
\phi$ as a 
composition of functions
$$\phi =z(1+\bar {b}a^{-1}z^{-n})^{1/n}\circ Rot(a^{-\frac 2n})\circ 
z(1+a^{-1}bz^n)^{-1/n}\tag 2.1.21$$
(This is the rigorous expression corresponding to the heuristic 
expression
$$\phi =(\frac {\bar {b}+\bar {a}z^n}{a+bz^n})^{1/n};\tag 2.1.22$$
one can also obtain (2.1.21) from (2.1.20) by using formulas from 
[FLM], especially Prop.  8.3.10, page 186.).  

In terms of the determinant formulas above, we have 
$$a^{\frac 2nL_0+\frac 1{12}n(n^2-1)\kappa}=detA_a(\phi )^{\kappa}
(\frac {detA_a(\phi )}{detA_p(\phi )})^{8L_0}\tag 2.1.23$$
This implies
$$det\vert A_a(\phi )\vert^2=(1-r^2)^{\frac 1{24}n(n^2-1)}$$
$$det\vert A_p(\phi )\vert^2=det\vert A_a(\phi )\vert^2(1-r^2)^{1
/4n}\tag 2.1.24$$
where $r=\vert a^{-1}b\vert$. As a consequence, if $\psi =\left(\matrix 
\alpha&\beta\\\bar{\beta}&\bar{\alpha}\endmatrix \right)$, then
$$det(A_a(\phi )A_a(\psi )A_a(\phi\circ\psi )^{-1})=(1+w\zeta )^{\frac 
1{12}n(n^2-1)},\tag 2.1.25$$
where $w=a^{-1}b$, $\zeta =\bar{\beta}\alpha^{-1}$.  It is a highly nontrivial matter to 
verify (2.1.24) directly.  

\smallskip

\flushpar Proof of (2.1.24) for $n=1$.  The antiperiodic determinant is 
1, since $\phi$ commutes with the polarization.  So we need only 
consider the periodic case.  

If we multiply $\phi$ on the left or right by a diagonal element, 
then $det\vert A_p\vert^2$ is unchanged.  Hence we can assume that $
\phi$ has 
the form $\phi =(1-r^2)^{-1/2}\left(\matrix 1&r\\r&1\endmatrix \right
)$, where $r<1$.  We have the 
Riemann-Hilbert factorization 
$$\Psi'(\theta )=\frac {\frac d{d\theta}(\frac {-r+e^{i\theta}}{1
-re^{i\theta}})}{i\frac {-r+e^{i\theta}}{1-re^{i\theta}}}=\frac {
1-r^2}{1+r^2-2rcos(\theta )}=(1-rz^{-1})^{-1}(1-r^2)(1-rz)^{-1}\tag 2.1.26$$
where $z=e^{i\theta}$.  Now $U_{\phi}(\cdot )=M_{\sqrt {\Psi'}}\Cal C_{
\psi}$, where $\Cal C_{\psi}f=f\circ\psi$.  The 
operator $\Cal C_{\psi}$ is nearly diagonal with respect to the Hardy 
polarization of $H$; to be precise, $C(\Cal C_{\psi})=0$, and $B(
\Cal C_{\psi})$ has rank 
$\le 1$, and its range is contained in $\Bbb Cz^0$. Using this, we see that
$$\aligned
A(U_{\phi})^{*}A(U_{\phi})&=A(\Cal C_{\psi})^{*}A(\sqrt {\Psi'})^
2A(\Cal C_{\psi})\\&=A(\Cal C_{\phi}\Psi'{}^{-1})A(\sqrt {\Psi'})^
2A(\Cal C_{\psi})\\&=A(\Cal C_{\phi})A(\Psi^{\prime -1})A(\sqrt {
\Psi'})^2A(\Cal C_{\psi})+B(\Cal C_{\phi})C(\Psi^{\prime -1})A(\sqrt {
\Psi'})^2A(\Cal C_{\psi})\\&=A(\Cal C_{\phi})(1+B(\Cal C_{\phi})C
(\Psi^{\prime -1})A(\Psi^{\prime -1})^{-1})A(\Psi^{\prime -1})A(\sqrt {
\Psi'})^2A(\Cal C_{\phi})^{-1}.\endaligned
\tag 2.1.27$$
Hence
$$det\vert A(U_{\phi})\vert^2=det(1+B(\Cal C_{\phi})C(\Psi^{\prime 
-1})A(\Psi^{\prime -1})^{-1})det(A(\Psi^{\prime -1})A(\sqrt {\Psi'}
)^2).\tag 2.1.28$$

We now have two determinants to evaluate.  Using the 
Riemann-Hilbert factorization (2.1.26), we see that 
$$A(\Psi^{\prime -1})=A((1-rz^{-1})(1-r^2)^{-1}A((1-rz)),\tag 2.1.29$$
and there is a similar decomposition for $A(\sqrt {\Psi'})$.  By inserting 
and rearranging terms, it is easy to see that 
$$det(A(\sqrt {\Psi'})^2A(\Psi^{\prime -1}))=det(S_1)^3\tag 2.1.30$$
where $S_1$ is the group-theoretic commutator
$$S_1=A((1-rz^{-1})^{-1/2})A((1-rz)^{-1/2})A((1-rz^{-1})^{1/2})A(
(1-rz)^{1/2}).\tag 2.1.31$$
The Helton-Howe formula for the determinant of such a 
commutator implies that (2.1.30) equals 
$$exp(3tr[log((1-rz^{-1})^{-1/2}),log((1-rz)^{-1/2}])=det(S_2)^{3
/4},\tag 2.1.32$$
where $S_2$ is the group-theoretic commutator obtained by 
removing all four square roots in the expression (2.1.31) for $S_
1$.
We now directly calculate that for $n>0$, $S_2z^n=z^n$, and
$$S_2z^0=proj_{_{H^{+}}}(z^0+(1-rz^{-1})^{-1}(1-rz)^{-1}rz^{-1})=\frac 
1{1-r^2}z^0+O(z).\tag 2.1.33$$
It follows that $det(S_2)=(1-r^2)^{-1}$, and hence (2.1.32) equals 
$(1-r^2)^{-3/4}$.  

For the other determinant, we calculate 
$$A(\Psi^{\prime -1})^{-1}z^0=\frac {1-r^2}{1-rz}proj_{H^{+}}(\frac 
1{1-rz^{-1}})=\frac {1-r^2}{1-rz}$$
$$C(\Psi^{\prime -1})A(\Psi^{\prime -1})^{-1}z^0=proj_{H^{-}}((1+
r^2-r(z+z^{-1}))\frac 1{1-rz})=-rz^{-1}$$
$$B(\Cal C_{\phi})(-rz^{-1})=(-r)proj_{H^{+}}(\frac {1+rz}{r+z})=
-r^2z^0.\tag 2.1.34$$
Since $Range(B(\Cal C_{\phi}))\subset \Bbb Cz^0$, it follows that
$$det(1+B(\Cal C_{\phi})C(\Psi^{\prime -1})A(\Psi^{\prime -1})^{-
1})=1-r^2.\tag 2.1.35$$
Thus (2.1.24) now follows from (2.1.28), the calculation of 
(2.1.30) in the preceding paragraph, and (2.1.35), in case $n=1$.  
// 

\bigskip

\centerline{\S 3. The measures $\nu_{\beta ,c,h}$, $\beta >0,c,h\ge 
0$.}

\bigskip

We begin by recalling the Malliavin-Shavgulidze construction of 
probability measures $\nu_{\beta}$ on $\Cal D_{C^1}$, indexed by inverse temperature 
$\beta$, which are quasi-invariant with respect to the left action of 
$\Cal D_{C^3}$.  The family of measures we consider differs slightly 
from the family considered in [K] and [MM], the point being 
that the family we consider has some chance of being 
asymptotically invariant as $\beta\downarrow 0$.  The quasi-invariance and 
expression for the Radon-Nikodym derivative follow from the 
general theory of Gaussian measures, which we recall in the 
first subsection.  

\bigskip

\flushpar\S 3.1. Preliminaries on Gaussian measures.

\smallskip

Let $H\to B$ denote an abstract Wiener space, and $\nu_G$ the 
associated Gaussian measure.  It is well-known that the affine 
automorphisms of $H$ (the Cameron-Martin space for $\nu_G$) fix the 
measure class of $\nu_G$, and the family of probability measures 
$$d\nu_{G,\beta}(b)=d\nu_G(\sqrt {\beta}b)=``\frac 1{\Cal Z}e^{-\frac {
\beta}2\langle b,b\rangle}\Cal Db\text{''}\tag 3.1.1$$
is asymptotically invariant with respect to these 
automorphisms in the following precise sense:  given an affine 
automorphism $\phi\in O(H)\propto H$, 
$$\int\vert\frac {d\phi_{*}\nu_{G,\beta}}{d\nu_{G,\beta}}-1\vert^
pd\nu_{G,\beta}\le 2\Gamma (\frac {p+1}2)(\beta\vert h_0\vert^2_H
)^{p/2}\tag 3.1.2$$
where $h_0$ is the translational part of $\phi$ (see (4.1.3) of Part III 
of [Pi]).

We need a nonlinear version of these results.  Part of what 
we need is standard, namely we have the following result on 
nonlinear transformations of Wiener space:  suppose that 

\smallskip

(a) $K:B\to H$ is a $C^1$ map with the property that the restriction 
of the derivative
$$dK\vert_b:H\to H\tag 3.1.3$$
is Hilbert-Schmidt, for each $b\in B$, and

(b) $T=1+K:B\to B$ is a homeomorphism, and the restriction of 
the derivative
$$dT\vert_b=1+dK\vert_b:H\to H\tag 3.1.4$$
is in $GL(H)$, for each $b\in B$; 

\smallskip

\flushpar then $T_{*}\nu_{G,\beta}$ is equivalent to $\nu_{G,\beta}$.  From the heuristic 
expression (3.1.1) one calculates that the Radon-Nikodym derivative 
is given heuristically by the formula 
$$\frac {d\nu_{G,\beta}(Tb)}{d\nu_{G,\beta}(b)}=det\vert dT\vert_
b\vert exp(-\frac {\beta}2(2\langle Kb,b\rangle -\vert Kb\vert_H^
2));\tag 3.1.5$$
the correct mathematical expression has the form
$$det_2\vert dT\vert_b\vert exp(-(``\beta\langle Kb,b\rangle -tr(
dK\vert_b)\text{''}-\frac {\beta}2\vert Kb\vert^2_H)),\tag 3.1.6$$
where the expression in parentheses is interpreted 
stochastically, and $det_2$ refers to the regularized 
Hilbert-Schmidt determinant for $dT\vert_b$ restricted to $H$ (see [R], 
especially section 4; in our application ordinary determinants 
suffice, and in this context the original result is due to Gross).  

It would be interesting to determine necessary and 
sufficient conditions under which 
$$\int\vert\frac {d\nu_{G,\beta}(Tb)}{d\nu_{G,\beta}(b)}-1\vert^p
d\nu_{G,\beta}(b)\to 0\quad as\quad\beta\downarrow 0,\tag 3.1.7$$
for each finite $p$.  It's very plausible that $spec(dK\vert_b)=\{
0\},$ 
$a.e.b$, is a necessary condition, but whether there is a useful 
general sufficient condition is unclear to me. We will consider 
this question in our specific context.

\bigskip

\flushpar\S 3.2. The Malliavin-Shavgulidze measure $\nu_{\beta}$.

\smallskip

The definition of $\nu_{\beta}$ involves a number of identifications.  
First, for any reasonable smoothness condition, there is a 
bijection 
$$\Cal D_1\times Rot(S^1)\leftrightarrow \Cal D:\psi_1,Rot(\lambda 
)\leftrightarrow\psi\tag 3.2.1$$
where $\psi^{-1}(1)=\lambda^{-1}$, $\psi_1=\psi\circ Rot(\psi^{-1}
(1))$, and $\Cal D_1$ denotes the 
stabilizer subgroup at $1\in S^1$.   This bijection induces an 
identification of $\Cal D_1$ with the left-coset space $\Cal D/Ro
t(S^1)$.  The 
left action is given explicitly by 
$$\Cal D\times \Cal D_1\to \Cal D_1:\phi ,\psi_1\to (\phi\circ\psi_
1)_1=\phi\circ\psi_1\circ Rot(\psi_1^{-1}(\phi^{-1}(1))).\tag 3.2.2$$
In turn there is an identification (where we now impose a 
specific smoothness condition) 
$$(\Cal D_{C^1})_1\leftrightarrow Path^{0,0}_{C^0}\Bbb R:\psi_1\leftrightarrow 
b=b_{\psi_1},\tag 3.2.3$$
where
$$b(t)=ln\Psi_1'(t)-ln\Psi_1'(0),\quad\Psi_1(t)=(\frac 1{2\pi}\int_
0^{2\pi}e^b)^{-1}\int_0^te^{b(\tau )}d\tau .\tag 3.2.4$$
Below we will routinely view $b$ as a $2\pi$-periodic function.  
Also it is useful to note that the map $\psi_1\to b_{\psi_1}$ extends to a 
map 
$$\Cal D_{C^1}\to Path^{0,0}_{C^0}\Bbb R:\psi\to b_{\psi}=ln\Psi'
-ln\Psi'(0);\tag 3.2.5$$
this extension satisfies the equation 
$$b_{\phi\circ\psi}=b_{\phi}\circ\Psi +b_{\psi}-(ln\Phi'(\Psi (0)
)-ln\Phi'(0)).\tag 3.2.6$$

In terms of the identification (3.2.3), the left action (3.2.2) 
becomes 
$$\Cal D_{C^1}\times Path^{0,0}_{C^0}\Bbb R\to Path^{0,0}_{C^0}\Bbb R
:\phi ,b\to b_{(\phi\circ\psi_1)_1}\tag 3.2.7$$
In particular for $\phi =Rot(\lambda ),$ $\lambda =e^{is}$,
$$\phi :b(t)\to b(t+T)-b(T),\tag 3.2.8$$
where $T=T_{\phi}(b)=\Psi^{-1}_1(-s)$, i.e. $T$ is the unique solution of
$$(\frac 1{2\pi}\int_0^{2\pi}e^b)^{-1}\int_0^Te^b=2\pi -s;\tag 3.2.9$$
and for $\phi\in (\Cal D_{C^1})_1$, 
$$\phi :b(t)\to b(t)+b_{\phi}\circ\Psi_1.\tag 3.2.10$$

Let $\nu_{\beta}^{0,0}$ denote the Brownian bridge probability measure on 
$B=Path^{0,0}_{C^0}\Bbb R$.  This is the Gaussian measure $\nu_{G
,\beta}$ corresponding 
to the Cameron-Martin Hilbert space $H=Path^{0,0}_{W^1}\Bbb R$, with 
$$\langle x,y\rangle_H=\int_0^{2\pi}x'(\tau )dy(\tau ).\tag 3.2.11$$
The following is essentially due to the Malliavins.  

\proclaim{(3.2.12)Proposition}The measure $\nu_{\beta}^{0,0}$ is quasi-invariant with 
respect to the left action of $\Cal D_{C^3}$ on $(\Cal D_{C^1})_1
\equiv Path^{0,0}_{C^0}\Bbb R$, and the 
Radon-Nikodym derivative is given by 
$$\frac {d\nu_{\beta}^{0,0}(b_{(\phi\circ\psi_1)_1})}{d\nu_{\beta}^{
0,0}(b)}=exp(-\frac {\beta}2\int_0^{2\pi}(b'_{\phi}(\Psi_1(\tau )
)^2-2b^{\prime\prime}_{\phi}(\Psi_1(\tau )))\{\frac {e^{b(\tau )}}{\frac 
1{2\pi}\int_0^{2\pi}e^b}\}^2d\tau ),$$
where $b=b_{\psi_1}$. 
\endproclaim

\flushpar(3.2.13)Remarks.  (a) Because $b_{Rot(\lambda )\circ\phi}
=b_{\phi}$ (by (3.2.6)), the 
formula above asserts that in particular $\nu_{\beta}^{0,0}$ is strictly 
invariant with respect to the left action (3.2.8) of $Rot(S^1)$ on the 
based loop space.  

(b) The analysis of Kosyak in [K] appears to correctly show 
that the induced representation 
$$(\Cal D_{C^3})_1\times L^2((\Cal D_{C^1})_1,d\nu_{\beta}^{0,0})$$
is irreducible. 

\smallskip

\flushpar Proof of (3.2.12).  We first prove invariance with 
respect to a rotation $\phi =Rot(e^{is})$.  Given a constant $T_0$, the 
Brownian bridge is invariant with respect to the time 
translation 
$$(T_0)_{*}:B\to B:b(t)\to b(t+T_0)-b(T_0),$$
since this defines a unitary transformation of the 
Cameron-Martin subspace.  As an immediate consequence, if we 
have a function $T:B\to \Bbb R$ having a finite range, then the 
corresponding (random translation of time) transformation, 
$$T_{*}:B\to B:b(t)\to b(t+T(b))-b(T(b)),$$
is an isomorphism of measure spaces. 

Now consider the transformation (3.2.8) defined by $\phi$.  For each 
$n$ we can define an approximation to $T_{\phi}$ having a finite number 
of values, say 
$$T_n:B\to \Bbb R:b\to k/n,\quad where\quad 2\pi k/n\le T_{\phi}(
b)\le 2\pi (k+1)n.$$
Each $(T_n)_{*}$ is measure-preserving. Given a bounded functional $
\Phi$ 
of the Brownian bridge, we have
$$\int\Phi\circ (T_n)_{*}\to\int\Phi\circ (T_{\phi})_{*}\quad as\quad 
n\to\infty ,$$
by dominated convergence. From this it follows that $(T_{\phi})_{
*}$ is 
measure-preserving.

Now suppose that $\phi\in (\Cal D_{C^3})_1$.  Let $h=b_{\phi}$, and let $
K$ denote the 
transformation 
$$K:B\to H:b\to h\circ\Psi_1,\tag 3.2.14$$
where $b=b_{\psi}$. We calculate that
$$dK\vert_b(x)(t)=\frac d{d\epsilon}\vert_{\epsilon =0}h((\frac 1{
2\pi}\int_0^{2\pi}e^{b+\epsilon x})^{-1}\int_0^te^{b+\epsilon x})$$
$$=h'(\Psi_1(t))((\frac 1{2\pi}\int_0^{2\pi}e^b)^{-1}\int_0^te^bx
-(\frac 1{2\pi}\int_0^{2\pi}e^b)^{-2}(\frac 1{2\pi}\int_0^{2\pi}e^
bx)\int_0^te^b)$$
$$=\frac {h'(\Psi_1(t))}{\frac 1{2\pi}\int_0^{2\pi}e^b}(\int_0^te^
bx-\Psi_1(t)\frac 1{2\pi}\int_0^{2\pi}e^bx).\tag 3.2.15$$
From this it is obvious that
$$dK\vert_b:H\to H\tag 3.2.16$$
is a trace class operator.  Hence $\nu_{\beta}^{0,0}$ is quasi-invariant.

To compute the Radon-Nikodym derivative, we first show that 
the spectrum for the operator in (3.2.16) is $\{0\}$.  To see this, 
let $\lambda\ne 0$ and suppose that (3.2.15) equals $\lambda x(t)$.  It is then 
straightforward to check that $x$ is a solution of the first-order 
linear equation 
$$\rho^2e^b(x-c)=-\mu\rho'x+\mu\rho x',$$
where $c=(\int_0^{2\pi}e^bx)/(\int_0^{2\pi}e^b)$, $\rho =h'\circ\Psi$, $
\mu =\lambda\frac 1{2\pi}\int_0^{2\pi}e^b$. The 
general solution of () is given by
$$x(t)=(-\frac c{\lambda})\rho (t)e^{\lambda^{-1}h(\Psi (t))}\int_
0^te^{-\lambda^{-1}h(\Psi (\tau ))}e^{b(\tau )}d\tau +C\rho (t)e^{
\lambda^{-1}h(\Psi (t))}.\tag 3.2.17$$
The condition $x(0)=0$ implies that $C=0$.  The condition 
$x(2\pi )=0$ forces either $c=0$, in which case $x=0$, or $\rho (
2\pi )=0,$ 
i.e.  $h'(2\pi )=0$.  It follows that if $\phi$ is such that $h'(
2\pi )\ne 0$ 
(recall that $h=b_{\phi}$), then the spectrum is $\{0\}$.  But now it 
follows by the continuous dependence of the spectrum on $\phi$ 
that the spectrum is always zero.  

We can now explicitly compute the Radon-Nikodym derivative 
from the stochastic expression (3.1.6), 
$$exp(-\frac {\beta}2\int_0^{2\pi}\{2(Kb)'db+(Kb)^{\prime 2}d\tau 
\}),$$
$$=exp(-\frac {\beta}2\int_0^{2\pi}\{2h'(\Psi_b(\tau ))\frac {e^{
b(\tau )}}{\frac 1{2\pi}\int_0^{2\pi}e^b}db+(h'(\Psi_b(\tau ))\frac {
e^{b(\tau )}}{\frac 1{2\pi}\int_0^{2\pi}e^b})^2d\tau \}).\tag 3.2.18$$
The meaning of the stochastic integral is that one integrates 
by parts, and this gives the expression in (3.2.12).//

\smallskip

We now want to discuss the question of whether $\nu_{\beta}^{0,0}$ is 
asymptotically invariant, i.e.  we want to evaluate the limit as 
$\beta\downarrow 0$ of 
$$\int\vert\frac {d\nu_{\beta}^{0,0}(b_{(\phi\circ\psi_1)_1})}{d\nu_{
\beta}^{0,0}(b)}-1\vert^pd\nu_{\beta}^{0,0}(b)$$
$$=\int\vert exp(-\frac {\beta}2(\int_0^{2\pi}(h'(\Psi_{b/\sqrt {
\beta}})^2-2h^{\prime\prime}(\Psi_{b/\sqrt {\beta}}))\{\frac {e^{\frac 
b{\sqrt {\beta}}}}{\frac 1{2\pi}\int_0^{2\pi}e^{\frac b{\sqrt {\beta}}}}
\}^2)d\tau )-1\vert^pd\nu_1^{0,0}(b)\tag 3.2.19$$
This limit would be zero if we could answer in the affirmative 
the following relatively straightforward

\proclaim{(3.2.20)Question}For any positive integer $n$, do we 
have
$$\lim_{\beta\to 0}\int exp(n\beta^2\frac {\int e^{\frac {2b}{\beta}}}{
(\int e^{\frac b{\beta}})^2})d\nu^{0,0}_1(b)=1?$$
\endproclaim

Since $h^{\prime\prime}$ is bounded, an affirmative answer would allow us to 
apply dominated convergence to (3.2.19), and the limit would be 
zero. 

To make sense of (3.2.20) we first must show that the function
$$1\le\int exp(n\beta^2\frac {\int e^{2b/\beta}}{(\int e^{b/\beta}
)^2})d\nu^{0,0}_1(b)\tag 3.2.21$$
is actually finite.  Given finiteness, because the function is 
even, and unbounded as $\beta$ goes to infinity, it is very plausible 
that it is an increasing function of $\beta^2$, and has a minimum at 
$\beta =0$.  

Is it plausible that (3.2.21) is finite?  Suppose that $H$ is a 
positive Morse function.  Stationary phase implies that 
$$\int e^{-H/T}=\sum_{t:H(t)=minH}e^{-H(t)/T}(\frac {2\pi T}{H^{\prime
\prime}(t)})^{1/2}+O(T^{3/2}).\tag 3.2.22$$
If we formally apply this to $H=b-sup\{b\}$, with $T=\beta$, then 
we see that 
$$\frac {\int e^{2b/\beta}}{(\int e^{b/\beta})^2}=\frac {\int e^{
2(b-sup\{b\})/\beta}}{(\int e^{(b-sup\{b\})/\beta})^2}=\beta^{-1/
2}*const+O(\sqrt {\beta}),\tag 3.2.23$$
and hence (3.2.20) seems trivially true.  However Brownian 
paths behave very differently than Morse functions in the 
vicinity of their (unique) maxima.  In fact Pittman and Yor 
have proven that the random variables
$$\beta^{-2}\int e^{(b-sup\{b\})/\beta}\tag 3.2.24$$
converge in law as $\beta\downarrow 0$ ([PY]).  This makes it very unclear
whether (3.2.21) is actually finite.  

Now consider the principal bundle
$$\matrix \Cal D_{C^1}&\leftarrow&Rot(S^1)\\\downarrow\\\Cal D_{C^
1}/Rot(S^1)\endmatrix \tag 3.2.25$$
Since $\nu_{\beta}^{0,0}$, viewed as a measure on the base, is 
$Rot(S^1)$-invariant, there is a unique $Rot(S^1)$-biinvariant 
probability measure on $\Cal D_{C^1}$, $\nu_{\beta}$, which projects to $
\nu_{\beta}^{0,0}$.  

\proclaim{(3.2.26)Corollary}The measure $\nu_{\beta}$ is quasi-invariant 
with respect to the left action of $\Cal D_{C^3}$.  If the answer to 
(3.2.20) is affirmative, then we have 
$$\int\vert\frac {d\nu_{\beta}(\phi\circ\psi )}{d\nu_{\beta}(\psi 
)}-1\vert^pd\nu_{\beta}(\psi )\to 0\quad as\quad\beta\to 0,$$
for each finite $p$.
\endproclaim

\flushpar\S 3.3. Existence of the measures $\nu_{\beta ,c,h}$.

\smallskip

In this subsection our task is to make sense of the measure 
which can be written heuristically as
$$\aligned
d\nu_{\beta ;c,h}=\frac 1Edet\vert A_a(\phi )\vert^{2c}\vert diag
(\phi )\vert^{16h}d\nu_{\beta}(\phi )\endaligned
$$
$$=\frac 1Eexp(-ctr\vert C_a(\phi )\vert^2)det_2\vert A_a(\phi )\vert^{
2c}\vert diag(\phi )\vert^{16h}d\nu_{\beta}(\phi );$$
and rigorously (as will show below) as
$$\frac 1{E'}exp(-\frac c{8\pi^2}\lim_{\delta\downarrow 0}\iint_{
\delta <\vert t-s\vert}(K_{\phi}(s,t)^2-E_{\beta}K_{\phi}(s,t)^2)
)det_2\vert A_a(\phi )\vert^{2c}\vert diag(\phi )\vert^{16h}d\nu_{
\beta}\tag 3.3.1$$
where $E_{\beta}$ denotes expectation with respect to $\nu_{\beta}$, and the 
kernel $K_{\phi}$ is given by (1.3.9).  

We already know that $\vert diag(\phi )\vert$ is a well-defined function of 
$\phi\in \Cal D_{qs}$, and it is bounded by 1, by (1.5.17) and (1.5.23).  Thus 
providing that $h\ge 0$, 
$$\vert diag\vert^{16h}\tag 3.3.2$$
is a well-defined random variable, because $\nu_{\beta}$ is supported on 
$\Cal D_{C^{\alpha}}$, provided that $\alpha <3/2$, and this is certainly contained in 
$\Cal D_{qs}$. The Hilbert determinant is taken care of by the following

\proclaim{(3.3.3)Corollary of (1.3.5)} The function
$$\phi\to det_2\vert A_a(\phi )\vert^2$$
is well-defined, continuous and bounded by 1 on 
$\{\phi\in \Cal D_{C^1}:log\Phi'\in B^{1/4}\}$.  
\endproclaim

\flushpar Proof.  The function is well-defined and continuous by 
(1.3.5).  Because $\vert A\vert^2+\vert C\vert^2=1$, we have $0\le
\vert C\vert^2\le 1$.  It follows 
from this that 
$$0\le det_2\vert A\vert^2=det((1-\vert C\vert^2)e^{\vert C\vert^
2})\le 1,\tag 3.3.4$$
for the function $\lambda\to (1-\lambda )exp(\lambda )$ is bounded by 1 for $
\lambda\ge 0$. //

\smallskip

Thus the Hilbert determinant in (3.3.3) is a well-defined bounded 
random variable with respect to $\nu_{\beta}$.  Now consider 
$$\iint_{\vert s-t\vert >\delta}\vert K_{\phi}(s,t)\vert^2-E_{\beta}
\iint_{\vert s-t\vert >\delta}\vert K_{\phi}(s,t)\vert^2,\tag 3.3.5$$
This is a well-defined random variable with respect to $\nu_{\beta}$, for 
each $\delta >0$.

\proclaim{(3.3.6)Proposition} (3.3.5) has a limit in probability $
[\nu_{\beta}]$ 
as $\delta\to 0$.
\endproclaim

\flushpar(3.3.7)Remark.  The basic intuition is that, by Theorem 
(1.2.6), $\vert K_{\phi}\vert^2$ and $\vert\frac {b(t)-b(s)}{e^{i
t}-e^{is}}\vert^2$ have essentially the same singular 
behavior along the diagonal.  Hence the integral of one can be 
regularized if and only if the integral of the other can be 
regularized, by subtracting out the expectation.  If we write $b$ 
in terms of its sine series, 
$$b(t)=\sum_{n>0}b_nsin(\frac n2t),\tag 3.3.8$$
then $\nu^{0,0}_{\beta}$ corresponds to the product measure 
$$\prod_{n>0}\sqrt {\frac {\beta n^2}{2\pi}}e^{-\frac {\beta n^2}
2b_n^2}db_n,\tag 3.3.9$$
and
$$\lim_{\delta\downarrow 0}\iint_{\vert e^{is}-e^{it}\vert >\delta}
\vert\frac {b(t)-b(s)}{e^{is}-e^{it}}\vert^2=\sum_{n>0}n(b_n^2-E_{
\beta}(b_n^2)),\tag 3.3.10$$
where the latter sum is absolutely convergent in the $L^2$ sense.  
Moreover, by Kolmogorov's theorem, it converges for $a.e.$ $b$
(see \S 3.1 of Part III of [Pi]). 

We will reduce the proof of (3.3.6) to a linear calculation 
which is essentially equivalent to this.  

\smallskip

\flushpar Proof of (3.3.6).  Since $\beta$ does not play an important role, 
we fix it such that $E(b(s)b(t))=s(2\pi -t),$ $s\le t$, where $E=
E_{\beta}$.  It 
is convenient to introduce the following notation.  Let 
$A=(b(t)+b(s))/2$, $\Delta =t-s$, $I=\int_s^te^b$, $\alpha =(\frac 
1{2\pi}\int_0^{2\pi}e^b)^{-1}$, and 
$$F=F(s,t;b)=A-ln(\frac {sin(\frac 12\alpha I)}{\alpha sin(\frac 
12\Delta )}).\tag 3.3.11$$

By (1.4.8) we have 
$$\vert K_{\phi}\vert^2=\frac {\vert e^F-1\vert^2}{2-2cos\Delta}.$$
Therefore
$$(e^F-1)^2-F^2=\sum_{n\ge 3}\frac 1{n!}(2^n-2)F^n,\tag 3.3.12$$
$$\vert K_{\phi}\vert^2-\frac {F^2}{2-2cos\Delta}=\sum_{n\ge 3}\frac 
1{n!}\frac {(2^n-2)F^n}{2-2cos\Delta}.\tag 3.3.13$$

We now want to estimate the expected value of this difference.
Note that
$$ln\frac {sin(\frac {\alpha}2I)}{\alpha sin(\frac 12\Delta )}=b(
u)\tag 3.3.14$$
for some $u$ between $s$ and $t$.  We will use the following 
standard fact several times, 

\proclaim{(3.3.15)Lemma}There is a constant $c$ such that
$$E(\sup_{s\le u\le t}\vert\frac {b(s)+b(t)}2-b(u)\vert^n)\le c2^
n(n!!)(\Delta (2\pi -\Delta ))^{n/2}.$$
\endproclaim

\flushpar Proof of (3.3.15).  For a standard Brownian motion $x(t
)$, we 
have 
$$Prob\{\sup_{s\le u\le t}\vert x(s)-x(u)\vert >\lambda \}=\frac 
2{\sqrt {2\pi (t-s)}}\int_{\lambda}^{\infty}e^{-\frac 1{2(t-s)}y^
2}dy\tag 3.3.16$$
$$E\sup_{s\le u\le t}\vert x(s)-x(u)\vert^n=\int\lambda^nd\{\frac 
2{\sqrt {2\pi (t-s)}}\int_{\lambda}^{\infty}e^{-y^2/2(t-s)}dy\}$$
$$=\left(\matrix n!!(t-s)^{^{n/2}},&n\quad even\\\sqrt {2/\pi}n!!
(t-s)^{n/2},&n\quad odd\endmatrix \right)\le n!!\Delta^{n/2}.\tag 3.3.17$$
(see Thm 2.5 of [Kn]). Similarly
$$E\sup_{s\le u\le t}\vert x(t)-x(u)\vert^n\le n!!\Delta^{n/2}.\tag 3.3.18$$

 Since $b(t)=x(t)-tx(2\pi )$, we have
$$\sup_{s\le u\le t}\vert b(s)-b(u)\vert^n=2^n\sup_u\vert x(s)-x(
u)+(u-s)x(2\pi )\vert^n$$
$$=2^{n-1}(\sup_u\vert x(s)-x(u)\vert^n+\Delta^n\vert x(2\pi )\vert^
n).\tag 3.3.19$$
Taking the expected value and using (3.3.17), we obtain
$$E\sup_{s\le u\le t}\vert b(s)-b(u)\vert^n\le 2^{n-1}(n!!\Delta^{
n/2}+n!!\Delta^n)\tag 3.3.20$$
There is a similar estimate for $E\sup_{s\le u\le t}\vert b(t)-b(
u)\vert^n$ (using 
(3.3.18)).  We also get similar estimates in terms of $(2\pi -\Delta 
)$ 
(using the periodicity of $b$).  It follows that 
$$E\sup_{s\le u\le t}\vert\frac {b(s)+b(t)}2-b(u)\vert^n\le E\frac 
12(\sup_u\vert b(s)-b(u)\vert^n+\sup_u\vert b(t)-b(u)\vert^n)$$
$$\le c2^nn!!(\Delta (2\pi -\Delta ))^{n/2}.//\tag 3.3.21$$

\smallskip
 
It follows from (3.3.14) and (3.3.15) that 
$E\vert F\vert^n\le cn!!(\Delta (2\pi -\Delta ))^{n/2}$.  Hence 
$$E(\iint_{0\le s,t\le 2\pi}\vert\vert K_{\phi}\vert^2-\frac {F^2}{
2-2cos\Delta}\vert )\le\sum_{n\ge 3}\frac {2^n-2}{n!}E\iint\frac {
\vert F\vert^n}{2-2cos\Delta}$$
$$\le\sum_{n\ge 3}\frac {2^n-2}{n!}c(n-1)!!\iint\frac {(\Delta (2
\pi -\Delta ))^{n/2}}{2-2cos\Delta}<\infty .\tag 3.3.22$$

Thus to prove (a) it suffices to prove that 
$$\lim_{\delta\to 0}\iint_{\vert 1-e^{i\Delta}\vert >\delta}\frac {
F^2-E(F^2)}{2-2cos\Delta}\tag 3.3.23$$
exists.  In the remainder of the proof, we will split the region 
of integration into two pieces, $\{\Delta\le\pi \}$ and the complement.  We 
will write out the estimates for the case $\Delta\le\pi$ only, since 
they are essentially the same for the complement, using 
periodicity.  

We write 
$$F^2=(A-ln\frac {\frac 12I}{sin(\frac {\Delta}2)}-ln\frac {sin(\frac {
\alpha}2I)}{\frac {\alpha}2I})^2$$
$$=(A-ln\frac {\frac 12I}{sin(\frac {\Delta}2)})^2-2(A-ln\frac {\frac 
12I}{sin(\frac {\Delta}2)})ln\frac {sin(\frac {\alpha}2I)}{\frac {
\alpha}2I}+(ln\frac {sin(\frac {\alpha}2I)}{\frac {\alpha}2I})^2.\tag 3.3.24$$
By (3.3.15), and the fact that $I$ involves integrating (not 
averaging) over an interval of length $\Delta$, 
$$E\{(A-ln\frac {\frac 12I}{sin(\frac {\Delta}2)})^2\}\le c\Delta 
,\quad E\{(ln\frac {sin(\frac {\alpha}2I)}{\frac {\alpha}2I})^2\}
\le c\Delta^4.\tag 3.3.25$$
Hence using (3.3.14), we have 
$$E\iint_{\Delta\le\pi}\vert\frac {\vert A-ln\frac {I/2}{sin(\Delta 
/2)}\vert^2-\vert A-ln\frac {sin\frac {\alpha}2I}{\alpha sin(\Delta 
/2)}\vert^2}{2-2cos\Delta}\vert <\infty .$$
Using similar estimates to replace $sin(\Delta /2)$ by $\Delta /2$, we see that 
it suffices to show that 

$$\lim_{\delta\downarrow 0}\iint_{\vert 1-e^{i\Delta}\vert >\delta 
,\Delta\le\pi}\frac {\vert A-ln(\frac 1{\Delta}I)\vert^2-E\vert A
-ln(\frac 1{\Delta}I)\vert^2}{2-2cos\Delta}dsdt\tag 3.3.26$$
exists for $a.e.$ $b$. 

Now we can write
$$A-ln(\frac 1{\Delta}I)=ln(\frac 1{\Delta}\int_s^te^{A-b(\tau )}
d\tau )=\frac 1{\Delta}\int_s^t(A-b)+ln(f),\quad f=\frac {\frac 1{
\Delta}\int_s^te^{A-b}}{e^{\frac 1{\Delta}\int_s^t(A-b)}},\tag 3.3.27$$
$$(A-ln(\frac 1{\Delta}I))^2=(\frac 1{\Delta}\int_s^t(A-b))^2+2(\frac 
1{\Delta}\int_s^t(A-b))ln(f)+ln^2(f).\tag 3.3.28$$
The important points are that the first term in (3.3.28) is 
linear in $b$, and the nonlinear term $f$ is now relatively 
tractable.  The following lemma states that we can regularize 
the first term in (3.3.28); it is essentially equivalent to (3.3.10).  

\proclaim{(3.3.29)Lemma}The
$$\lim_{\delta\downarrow 0}\iint_{\delta <\Delta <\pi}\frac {\vert\frac 
1{\Delta}\int_s^t(A-b)\vert^2-E\vert\frac 1{\Delta}\int_s^t(A-b)\vert^
2}{\Delta^2}$$
exists in the $L^2$ sense. 
\endproclaim

Assuming this for the moment, we can complete the proof of 
(3.3.6) by showing that the last two terms in (3.3.28) do not 
require any regularization, i.e.  
$$\iint_{\Delta <\pi}\frac {E((\frac 1{\Delta}\int_s^t\vert A-b\vert 
)ln(f)+ln^2(f))}{\Delta^2}<\infty .\tag 3.3.30$$

We have
$$E((\frac 1{\Delta}\int_s^t(A-b))^2)\le\frac 1{\Delta}\int_s^tE(
A-b)^2\le\Delta .\tag 3.3.31$$
Because $f\ge 1$, $ln^2f\le (f-1)^2$, hence
$$E(ln^2(f))\le E((f-1)^2)=1+Ef^2-2Ef.\tag 3.3.32$$
Now
$$f=(1+\frac 1{\Delta}\int_s^t(A-b)+\frac 12\frac 1{\Delta}\int_s^
t(A-b)^2+..)(1-\frac 1{\Delta}\int_s^t(A-b)+\frac 12(\frac 1{\Delta}
\int_s^t(A-b))^2-..)\tag 3.3.33$$
Using
$$E(\frac 1{\Delta}\int_s^t(A-b))^n\le E\frac 1{\Delta}\int_s^t(A
-b)^n=\frac 1{\Delta}\int_s^tE(A-b)^n\le c\Delta^{n/2},\tag 3.3.34$$
we see that
$$Ef=1+\frac 12E(\frac 1{\Delta}\int_s^t(A-b)^2-(\frac 1{\Delta}\int_
s^t(A-b))^2)+O(\Delta^{3/2}).\tag 3.3.35$$
Similarly,
$$Ef^2=1+E(\frac 1{\Delta}\int_s^t(A-b)^2-(\frac 1{\Delta}\int_s^
t(A-b))^2)+O(\Delta^{3/2}).\tag 3.3.36$$
It follows that (3.3.32) is $O(\Delta^{3/2})$.  This together with (3.3.31) 
implies (3.3.30).//  

\smallskip

\flushpar Proof of (3.3.29).  Since this is a linear problem, we 
will just indicate what is involved.  We must show that in the 
$L^2$ sense, 
$$\iint_{\delta <\Delta <\epsilon}\frac {\vert\frac 1{\Delta}\int_
s^t(A-b)\vert^2-E\vert\frac 1{\Delta}\int_s^t(A-b)\vert^2}{\Delta^
2}\to 0\quad as\quad\delta ,\epsilon\to 0,\quad i.e.$$
$$E((\iint_{\delta <\Delta <\epsilon}(\frac 1{\Delta^2}\int_s^t(A
-b))^2)^2)-(E\iint_{\delta <\Delta <\epsilon}\vert\frac 1{\Delta^
2}\int_s^t(A-b)^2)^2\to 0\quad in\quad L^2.$$
By expanding all the terms, all the expectations involved here 
can be computed explicitly.//  

\smallskip

We can now formulate our main

\proclaim{(3.3.37)Theorem}For $c,h\ge 0$, (a) there exists an essentially 
unique measure $\nu_{\beta ;c,h}$ which (i) is absolutely continuous with 
respect to $\nu_{\beta}$, and (ii) has Radon-Nikodym derivative 
$$\frac {d\nu_{\beta ;c,h}(\phi\circ\psi )}{d\nu_{\beta ;c,h}(\psi 
)}=\frac {d\nu_{\beta}(\phi\circ\psi )}{d\nu_{\beta}(\psi )}\vert\frac {
\sigma_{c,h}(\hat{\phi}\circ\hat{\psi })}{\sigma_{c,h}(\hat{\psi }
)}\vert^2$$
for all $\phi\in \Cal D_{C^3}$; and (b) the measure $\nu_{\beta ;
c,h}$ is finite, so that we 
can normalize it to be a probability measure. 

\endproclaim

\flushpar(3.3.38)Remark. Note that
$$\vert\frac {\sigma_{c,h}(\hat{\phi}\circ\hat{\psi })}{\sigma_{c
,h}(\hat{\psi })}\vert^2=\vert\frac {diag(\phi\circ\psi )}{diag(\psi 
)}\vert^{16h}det\vert A_a(\phi )+B_a(\phi )Z_a(\psi )\vert^2,$$
where $Z_a=C_aA_a^{-1}$. With probability one $Z_a(\psi )\in \Cal L_
p$ for all $p>2$, 
by (1.3.5), hence the produce $B_a(\phi )Z_a(\psi )$ is trace class, and the 
determinant is well-defined.

\smallskip

\flushpar Proof of (3.3.37).  Let $\nu_{\beta ;c,h}$ be the measure given by 
the last line of (3.3.3), where we interpret the exponential 
term as a random variable using (3.3.6), and we temporarily set 
$E'=1$.  It is a simple matter to check that it has the 
Radon-Nikodym derivative claimed above.  It is essentially 
unique because $\nu_{\beta}$ is ergodic with respect to the left action of 
$\Cal D_{C^3}$, as established by Kosyak in [K]. //

\smallskip

\flushpar(3.3.39)Remark.  It is easy to check that $\nu_{\beta}$ is right 
quasi-invariant only with respect to rotations, with respect to 
which it is strictly invariant.  Hence the same is true for 
$\nu_{\beta ;c,h}$.  Thus these regularized Virasoro measures, compared to 
those for loop groups, are fundamentally asymmetric.  

\bigskip

\centerline{\S 4. On Existence of Invariant Measures}

\bigskip

Via the natural inclusion
$$\Cal D_{qs}\to \Cal N^{-}\cdot \Cal H\cdot \Cal N^{+}\tag 4.1.1$$
given by (1.5.17), we can view each measure $\nu_{\beta ;c,h}$ as a probability 
measure on $\Cal N^{-}\cdot \Cal H\cdot \Cal N^{+}$. 

Given $\phi\in \Cal D_{qs}$, we write
$$\phi =l\circ diag\circ u,\tag 4.1.2$$
as in (1.5.17), and
$$u(z)=z(1+\sum_{n\ge 1}u_nz^n),\quad l^{-1}(z)=z+\sum_{n\ge 0}b_
nz^{-n}.\tag 4.1.3$$
By (1.5.23) and the Bieberbach-de Branges inequalities, we know 
that $\vert diag\vert\le 1$, $\vert b_n\vert\le 1/n$, and $\vert 
u_n\vert\le n+1$.  

\proclaim{(4.1.4)Corollary}For each $c,h\ge 0$, the family of measures
$$\{\nu_{\beta ,c,h}:\beta >0\}$$
has weak limits as $\beta\to 0$, with respect to $BC(\Cal N^{-}\cdot 
\{\vert z\vert\le 1\}\cdot \Cal N^{+})$ 
(bounded continuous functions). Any such limit is supported on 
the space 
$$\prod \{\vert b_n\vert\le 1/n\}\cdot \{\vert diag\vert\le 1\}\cdot
\prod \{\vert u_n\vert\le n+1\}.$$
\endproclaim

The question now is whether the limit points are 
$\Cal D_{an}$-invariant, viewed as bundle-valued measures.  As it stands, 
this question is not well-posed, because $\Cal D_{an}$ does not act on 
the formal completion.  One way to proceed would be to first 
investigate whether the limit points are infinitesimally 
invariant.  However, ultimately one would like to know 
whether the limit points are supported on $\Cal D_{qs}$, so that the 
global question would be well-posed.  

To relate this to a concrete problem, it is convenient to
again refer to the case of loop groups.  In that case the 
analogue of the factorization (4.1.2) is the Riemann-Hilbert 
factorization 
$$g=g_{-}\cdot g_0\cdot g_{+}\tag 4.1.5$$
where $g_0\in G$, $g_{\pm}\in G(\Cal O(D^{\pm}))_1$, respectively. If we write
$$g_{+}=1+\hat {g}_{+}(1)z+...,\tag 4.1.6$$
then it is easy to prove that
$$\nu^{*=*}_{\beta}\{\vert\hat {g}_{+}(n)\vert >R\}\to 0\quad as\quad 
n\to\infty ,\tag 4.1.7$$
for each fixed $R$.  If we could prove that this limit is 
monotone (which, at least to me, seems intuitively reasonable), 
then we could describe the support of $\mu_0=\lim_{\beta\to\infty}
\nu_{\beta}^{*=*}$ quite 
precisely.  This depends upon the fact that we can estimate 
the distribution for the first coefficient $\hat {g}_{+}(1)$, using invariance 
considerations.  

For this method to work in the Virasoro case, we would need an 
estimate for the $u_1$ distribution (which we suspect is 
absolutely continuous only for $c>1$). Assuming that this has 
been resolved in some way, we would then need something like 
the following

\proclaim{(4.1.8)Conjecture}For each $\beta >0$, $c\ge 0$,
$$\nu_{\beta ,c}\{\vert c_n\vert >(n+1)R\}\le\nu_{\beta ,c}\{\vert 
c_{n-1}\vert >nR\}.$$
\endproclaim

This would clearly imply the Bieberbach-de Branges 
inequalities, since we know that $\vert u_1\vert <2$ with probability one.
It would follow from this conjecture that we could estimate the 
distribution of $u_n$ in terms of the distribution for $u_1$, and 
hopefully this would show that the limit points are supported 
on a space such as $\Cal D_{qs}$. 

Concerning the diagonal distribution, as in the case of loop 
groups, it should be possible to a priori arrive at a reasonable 
conjecture for the evaluation of the integral 
$$\int\dot {a}(\phi )^{-i\lambda}d\mu_{c,h}(\phi )=\lim_{\beta\downarrow 
0}\int\dot {a}(\phi )^{-i\lambda}d\nu_{\beta ,c,h}(\phi ),\tag 4.1.9$$
where we are writing
$$\phi =l\cdot diag\cdot u,\quad diag=m\dot {a},\quad 0<\dot {a}<
1,\tag 4.1.10$$
but I have not succeeded at this.

\bigskip

\flushpar Acknowledgements. I thank Nasser Towghi for many 
useful conversations, and Jan Wehr for pointing out [PY] to 
me.

\bigskip

\centerline{References}

\bigskip

\flushpar[B] L. de Branges, A proof of the Bieberbach 
conjecture, Acta Mathematica 154 (1985) 137-152.

\flushpar[FLM] I. Frenkel, J. Lepowsky, and A. Meurman, Vertex 
Operator Algebras and the Monster, Academic Press (1990).

\flushpar[K] A.V. Kosyak, Irreducible regular Gaussian 
representations of the group of the interval and circle 
diffeomorphism, J. Funct. Anal. 125 (1994) 493-547.

\flushpar[Kn] F. Knight, Essentials of Brownian Motion and 
Diffusion, Mathematical Surveys, No. 18, AMS (1981). 

\flushpar[KY] A.A.  Kirillov and D.V.  Yuriev, Representations 
of the Virasoro algebra by the orbit method, J.  Geom.  Phys.,  
Vol.  5, n.  3 (1988).  

\flushpar[L] O. Lehto, Univalent Functions and Teichmuller 
Spaces, Springer-Verlag (1986).

\flushpar[MM] M.  Malliavin and P.  Malliavin, Quasi-invariant 
measures on the group of diffeomorphisms of the circle, 
Proceedings of the Lisbonne Conference, ed.  by A.B.  Cruzeiro, 
Birkhauser (1991).  

\flushpar[MP] R. Moody and Pianzola, Lie algebras with Triangular 
Structure, Wiley (1995).

\flushpar[NS] S.  Nag and D.  Sullivan, Teichmuller theory and 
the universal period mapping via quantum calculus and the 
$H^{1/2}$ space on the circle, Osaka J.  Math.,  32, no.1 (1995).  

\flushpar[Pe] O. Pekonen, Universal Teichmuller space in 
geometry and physics, J. Geom. Phys. 15 (1995) 227-251.

\flushpar[Pi] D. Pickrell, Invariant measures for unitary forms 
of Kac-Moody Lie groups, Parts I-III, submitted to J. Funct. 
Anal., Funct-an 9510005

\flushpar[PY] J.W. Pitman and M. Yor, A limit theorem for 
one-dimensional Brownian motion near its maximum, and its 
relation to a representation of the 2-dimensional Bessel bridge, 
preprint.

\flushpar[PS] A. Pressley and G. Segal, Loop Groups, Oxford 
University Press (1986).

\flushpar[R] R. Ramer, On nonlinear transformations of Gaussian 
measures, J. Funct. Anal. 15 (1974) 166-187.

\flushpar[Se1] G. Segal, Unitary representations of some infinite 
dimensional groups, Comm. Math. Phys. 80 (1981), 301-342.

\flushpar[Se2] ----, The definition of conformal field theory, 
preprint.  

\flushpar[S] S. Semmes, On the Cauchy integral and related 
operators on smooth curves, dissertation, Washington 
University, St. Louis, Missouri (1983)

\flushpar[Sh] E.T. Shavgulidze, Distributions on infinite 
dimensional spaces and second quantization in string theories, 
II, in ``V International Vilnius Conference in Probability Theory 
and Mathematical Statistics'' (1989) 359-360.

\end